\documentclass[12pt]{article}
\usepackage{amsmath,
amssymb}
\newcommand{\matr}[4]{\left[ \begin{array}{cc}
                              #1 & #2 \\
                              #3 & #4
                              \end{array} \right]}
\newcommand{\ra}{\rightarrow}
\newcommand{\la}{\leftarrow}
\newcommand{\da}{\downarrow}
\newcommand{\lra}{\longrightarrow}
\newcommand{\lla}{\longleftarrow}
\newcommand{\sra}[1]{\stackrel{#1}{\ra}}
\newcommand{\slra}[1]{\stackrel{#1}{\lra}}
\newcommand{\sla}[1]{\stackrel{#1}{\la}}

\newcommand{\sbth}[3]{\substack{#1 \\ #2 \\ #3}}
\newcommand{\sea}{\searrow}
\newcommand{\swa}{\swarrow}
\newcommand{\ov}{\overline}

\newcommand{\Nm}{{\rm Nm}\,}
\newcommand{\Ima}{{\rm Im}\,}
\newcommand{\Ker}{{\rm Ker}\,}
\newcommand{\Prym}{{\rm Prym}\,}
\newcommand{\Jac}{{\rm Jac}\,}
\newcommand{\Pic}{{\rm Pic}\,}

\renewcommand{\Sp}{{\rm Sp}\,}
\newcommand{\snu}{\mbox{$\scriptstyle{\nu}$}}
\newcommand{\snutil}{\mbox{$\scriptstyle{\tilde{\nu}}$}}
\newcommand{\spi}{\mbox{$\scriptstyle{\pi}$}}
\newcommand{\spione}{\mbox{$\scriptstyle{\pi_1}$}}
\renewcommand{\sf}{\mbox{$\scriptstyle{f}$}}
\newcommand{\sg}{\mbox{$\scriptstyle{g}$}}
\newcommand{\sgone}{\mbox{$\scriptstyle{g_1}$}}
\newcommand{\sh}{\mbox{$\scriptstyle{h}$}}

\renewcommand{\deg}{{\rm deg}\,}

\newcommand{\bbB}{\mathbf}
\newcommand{\CC}{{\bbB C}}
\newcommand{\GG}{{\bbB G}}
\newcommand{\PP}{{\bbB P}}
\newcommand{\QQ}{{\bbB Q}}
\newcommand{\RR}{{\bbB R}}
\newcommand{\ZZ}{{\bbB Z}}
\newcommand{\Ztwo}{\ZZ/2\ZZ}
\newcommand{\cA}{{\cal A}}

\newcommand{\cM}{{\cal M}}
\newcommand{\cS}{{\cal S}}
\newcommand{\cT}{{\cal T}}
\newcommand{\cTb}{\ov{{\cal T}}}
\newcommand{\taub}{\ov{\tau}}

\newcommand{\Ctil}{\tilde{C}}
\newcommand{\Ttil}{\tilde{T}}
\newcommand{\Ptil}{\tilde{P}}
\newcommand{\nutil}{\tilde{\nu}}
\newcommand{\PrymC}{\Prym(\Ctil/C)}
\newcommand{\PrymnuC}{\Prym(\nu \Ctil/\nu C)}
\newcommand{\PrymCp}{\Prym(\tilde{{C'}}/C')}
\newcommand{\Ttilp}{\tilde{{T'}}}
\newcommand{\Tph}{T'_{h}}
\newcommand{\Tpt}{T'_{t}}
\newcommand{\Ttilph}{\tilde{\Tph}}
\newcommand{\Ttilpt}{\tilde{\Tpt}}
\newcommand{\Gm}{\GG_m}
\newcommand{\we}{{\rm w}}

\newcommand{\ssub}{\scriptstyle{\subset}}
\newcommand{\ssup}{\scriptstyle{\supset}}
\newcommand{\sslash}{\scriptstyle{/}}
\newcommand{\seq}{\scriptstyle{=}}
\newcommand{\smin}{\scriptstyle{-}}
\newcommand{\stimes}{\times}
\newcommand{\bto}{\ssub \seq \sslash \smin\smin }
\newcommand{\btt}{\stackrel{\subset}{\ssub} \sslash \ssub }
\newcommand{\btth}{\ssub \ssub \sslash \smin\smin }
\newcommand{\btf}{\ssup\!\ssub \sslash
               \stimes }
\newcommand{\tto}{\stackrel{\subset}{\ssub}\seq 
                             \sslash\ssub \smin }
\newcommand{\ttt}{\ssub\ssub\seq \sslash \smin\smin\smin }
\newcommand{\tttp}{\ssup\!\ssub\seq
            \sslash     \stimes\smin  }
\newcommand{\ttoX}{\ssub \smin\smin }
\newcommand{\tttX}{\ssub \ssub }

\newtheorem{lemma}{Lemma}
\newtheorem{proposition}[lemma]{Proposition}
\newtheorem{theorem}[lemma]{Theorem}

\newtheorem{remark}[lemma]{Remark}

\newenvironment{pf}{\noindent {\it Proof:}\/ }
                       {\par \medskip \par}

\begin{document}
\begin{titlepage}
\begin{center}
\hfill {\tt alg-geom/9712027}\\

\vskip .4in  

{\bf The arithmetic-geometric mean}\\
{\bf and isogenies for curves of higher genus}

\vskip .3in

Ron Donagi\footnote{Partially supported by NSF grant DMS 95-03249}

\vskip .1in

{\em Department of mathematics, University of Pennsylvania, Philadelphia, PA 19104-6395}

\vskip .1in

and

\vskip .1in

Ron Livn\'e

{\em Institute of mathematics, Hebrew university of Jerusalem, Givat Ram 91904, Israel }

\end{center}

\vskip .2in

\noindent
{\bf Abstract}\\

\noindent
Computation of Gauss's arithmetic-geometric mean involves iteration of a 
simple step, whose algebro-geometric interpretation is the construction of 
an elliptic curve isogenous to a given one, specifically one whose period 
is double the original period. A higher genus analogue should involve the 
explicit construction of a curve whose jacobian is isogenous to the jacobian 
of a given curve. The doubling of the period matrix means that the kernel of 
the isogeny should be a lagrangian subgroup of the group of points of order
$2$ in the jacobian. In  genus $2$ such a construction was given classically 
by Humbert  and was studied more recently by Bost and Mestre. In this 
article we give such a construction for general curves of  genus $3$. We 
also give a similar but simpler construction for hyperelliptic curves of genus 
$3$. We show that the hyperelliptic construction is  a degeneration of the 
general one, and we prove that the kernel of the induced isogeny on 
jacobians is a lagrangian subgroup of the points of order $2$. We show 
that for $g \geq 4$ no similar construction exists, and we also reinterpret the
genus $2$ case in our setup. Our construction of these correspondences 
uses the bigonal and the trigonal constructions, familiar in the theory of 
Prym varieties.

\end{titlepage}
\vfill
\eject


\section{Introduction}
 It is well-known that computation of Gauss's
arithmetic-geometric mean involves iteration of a simple
step, whose algebro-geometric interpretation is the
construction of an elliptic curve isogenous to a given
one, specifically one whose period is double the original
period (for a modern survey see \cite{cox}). A higher
genus analogue should involve the explicit
construction of a curve whose jacobian is isogenous to
the jacobian of a given curve. The doubling of the
period matrix means that the kernel of the isogeny should
be a lagrangian subgroup of the group of points of order
$2$ in the jacobian. In  genus $2$ such a construction was
given classically by Humbert \cite{hum} and was studied
more recently by Bost and Mestre \cite{bome}. In this article 
we give such a construction for general curves of 
genus $3$. We also give a similar but simpler 
construction for hyperelliptic curves of genus $3$.
We show that the hyperelliptic construction is 
a degeneration of the general one, and we prove that the kernel
of the induced isogeny on jacobians is a lagrangian subgroup
of the points of order $2$. We show that for $g \geq 4$ no similar
construction exists, and we also reinterpret the
genus $2$ case in our setup.

To construct these correspondences we use the bigonal
and the trigonal constructions, familiar in the
theory of Prym varieties (\cite{don}). In genus $2$
Bost and Mestre note that Humbert's construction induces 
on jacobians an isogeny whose kernel is of type
$(\Ztwo)^2$. We show that Humbert's construction
is an instance of the bigonal construction, and prove
that the above kernel is a lagrangian subgroup of the 
points of order $2$. In fact Bost and Mestre use Humbert's
construction to give a variant of Richelot's genus
$2$ arithmetic-geometric mean.
In light of the clear analogy, in particular the
fact that a generic principally polarized abelian
threefold is a jacobian, one might hope that our
construction could be used in a similar way.

We work throughout over an algebraically closed field of
characteristic $0$. However, our methods clearly extend more 
generally. For example, the results of Section~\ref{genustwo} 
hold if the characteristic is not $2$, and those of 
Sections~\ref{hyperelliptic} and \ref{genusthree} if it is $>3$.

The first author thanks the Hebrew University of Jerusalem
and the Institute for Advanced Studies in Princeton for their
hospitality during the time this work was done. 
The second author thanks the University of Pennsylvania for 
its hospitality while this article was being written.

\section{Preliminaries} \label{prel}
\noindent {\bf Polarizations.} For an abelian variety
$A$ denote by $A[n]$ the kernel of multiplication by $n$.
In the sequel we will need the following standard facts and
notation.
\begin{enumerate}
\item A polarization $\Theta$ on an abelian variety $A$
induces by restriction a polarization $\Theta_B$ on any
abelian subvariety $B$ of $A$.
\item Recall that the type of a polarization $\Theta$ on a
$g$\/-dimensional abelian variety is a g\/-tuple of
positive integers $d_g|\dots|d_2|d_1$. We say that
$\Theta$ is a principal polarization if it is of type
$1^g = (1,\dots,1,1)$ ($g$\/ times). In that case, suppose
that $p$ is a prime, and that $K$ is a subgroup of $A[p]$
isomorphic to $(\ZZ/p\ZZ)^r$ and isotropic for the Weil
pairing $\we_p$. Then $\Theta$ induces a polarization
on $A/K$, characterized by the property that its pull back
to $A$ is $p\Theta$. Its type is then $p^{g-r}\cdot 1^r$.
In this situation we will say that $K$ is a lagrangian
subgroup of $A[p]$ if $r=g$.
\item The type of a polarization is preserved under continuous
deformations.
\end{enumerate}

\vspace*{0.1cm}

\noindent {\bf Double covers.} Given a double cover, i.e. a
finite morphism $\pi: \Ctil \ra C$ of degree $2$ between smooth 
projective curves, the Prym variety $\Prym(\Ctil/C)$ is defined
to be the connected component of the kernel of the norm map
\[ \pi_*: \Jac(\Ctil) \ra \Jac(C). \]
It is an abelian variety, and it has a natural principal
polarization when $\pi$ is unramified, namely one half of
the polarization induced on it as an abelian subvariety of
$\Jac(\Ctil)$ (\cite{mum2}). This definition extends to 
singular curves $C$, $\Ctil$, if we interpret $\Jac$ as
the (not necessarily compact) generalized jacobian. This was
studied by Beauville \cite{bea}.
Particularly important for us will be the cases when
1. $C$, $\Ctil$ have only ordinary double points,
2. $\pi^{-1}(C_{\rm sing}) = \Ctil_{\rm sing}$, and
3. for each $x\in C_{\rm sing}$ the inverse
image $\pi^{-1}(x)$ consists of a single point,
and each branch of $\pi^{-1}(x)$ maps to a different branch
of $x$ and is ramified over it.
(We shall then say that $\pi$ is of Beauville type at $x$.)
In such cases $\Prym(\Ctil/C)$ is compact, and the following
three conditions are equivalent:
\begin{enumerate}
\item $\pi$ is unramified away from $C_{\rm sing}$.
\item The arithmetic genera satisfy $g(\Ctil) = 2g(C)-1$.
\item The cover $\Ctil/C$ is a flat limit of smooth unramified
double covers.
\end{enumerate}
We shall call a cover satisfying these conditions allowable;
from the third condition we see that the $\Prym$ is prinicipally
polarized in such a case.

Let $C$ be a curve having only ordinary double points as
singularities, and let $\nu_x: N_x \ra C$ be the normalization
map of exactly one such singular point $x$. We denote by $L(x)$
the line bundle of order $2$ in $\Ker\nu_x^\ast$. (It is
obtained from the trivial line bundle on $N_x$ by gluing the
fibers over the two inverse images of $x$ with a twist of $-1$
relative to the natural identification.)

\begin{lemma} \label{polar}
Let $\pi: \Ctil \ra C$ be an allowable double cover, $\nu \pi:
\nu \Ctil \ra \nu C$ its (partial) normalization at $r \geq 1$
ordinary double points $x_1,\dots,x_r$. Let $g$ be the (arithmetic)
genus of the partial normalization $\nu C$, so the
arithmetic genus of $C$ is $g+r$. Then $\Prym (\Ctil/C)$
has a principal polarization,
$\Prym(\nu \Ctil/ \nu C)$ has a polarization of type
$2^{g}1^{r-1}$, and the pullback map
$\nu^*:\Prym(\Ctil/C) \ra \Prym(\nu\Ctil/\nu C)$ is an isogeny
of degree $2^{r-1}$. The kernel of $\nu^*$ is the subgroup
of $\Prym(\Ctil/C)[2]$ generated by the pairwise differences
of the line bundles $L(x_i)$ defined above. This subgroup is
isotropic for the mod $2$ Weil pairing $w_2$.
\end{lemma}
\begin{pf} The generalized jacobians fit in short exact
sequences
\[\begin{array}{ccccccccc}
0 & \ra & \Gm^r & \ra & \Jac(\Ctil) & \ra & \Jac(\nu \Ctil) &\ra & 0 \\
  &     &  \da &     & \da         &     & \da             &    &    \\
0 & \ra & \Gm^r & \ra & \Jac(C)     & \ra & \Jac(\nu C)     &\ra & 0
\end{array} \]
where the vertical maps are the norm maps induced by $\pi$
and by $\nu\pi$. We
compare the kernels: to begin with, the kernel of the norm
map is connected for ramified double covers (in particular
for $\nu \Ctil/\nu C$), and has two components for unramified
covers. This is shown in \cite{mum2} in the nonsingular case,
and so by continuity this holds also for allowable singular
covers (in particular for $\Ctil/C$). The multiplicative
groups parametrize extension data
and the norm is the squaring map. So the short exact sequence
of kernels gives
\[0 \ra (\Ztwo)^r \ra \PrymC \times \Ztwo \stackrel{\nu^*}{\ra}
\PrymnuC \ra 0,\]
and the first part of the lemma follows.

To prove that the subgroup $(\Ztwo)^{r-1}$ of
$\PrymC[2]$ is isotropic for $\we_2$, notice that
its generators are reductions modulo $2$ of the
vanishing cycles for $\Ctil$, and vanishing cycles for distinct
ordinary double points are disjoint. Therefore these vanishing
cycles have $0$ intersection number in $\ZZ$\/- (or $\QQ$\/-)
homology. By the definition of the polarization of $\PrymC$
in Section~\ref{prel}, and the well-known expression for the
Weil pairing in terms of the intersection (or cup product)
pairing (see e.g. \cite[theorem 1, Ch. 23]{mum}), the
rest of the lemma follows.
\end{pf}

\section{The bigonal and the trigonal constructions} 
\label{bigonal}
There are several elementary constructions which associate a
double cover of some special kind with another cover (or curve)
with related Prym (of Jacobian).
We now review the bigonal and the trigonal constructions,
following (\cite{don}). Assume we
are given smooth projective curves $\Ctil$, $C$ and $K$ and
surjective maps $f:C\ra K$ and $\pi:\Ctil\ra C$, so that
$\deg\pi = \deg f=2$ over any component. The bigonal
construction associates new curves and maps of the same
type $\tilde{{C'}}\sra{\pi'} C'\sra{f'}K$  as follows. Let
$U\subset K$ be the maximal open subset over which $f\pi$
is unramified. Then $\tilde{{C'}}$ represents over $U$ the sheaf of
sections, in the complex or the \'etale topology, of
$\pi:(f\pi)^{-1}U\ra f^{-1}U$. It is a $4$\/-sheeted
cover of $U$. We then view $\tilde{{C'}}_{|U}$ as a locally
closed subvariety of $\Ctil \times \Ctil$ and define
$\tilde{{C'}}$ as the closure.
The projection to $U$ extends to a morphism $\tilde{{C'}}\ra K$,
and the involution $\iota$ of $\tilde{{C'}}_{|U}$
which sends a section to the complementary section extends
to $\tilde{{C'}}$. We define $C'=\tilde{{C'}}/\iota$ and $f'$ and
$\pi'$ as the quotient maps.

We will need to extend this construction  to allowable covers
of curves with ordinary double points; however in a family
acquiring a singularity of Beauville type the arithmetic genus
of the resulting $\tilde{{C'}}$ is not locally constant.
More technically, the naive construction as the
closure of $\tilde{{C'}}_{|U}$ in $\Ctil\times\Ctil$ is not
flat in families, which is not adequate for our purposes: for
example, we want the bigonal construction to be symmetric.

To achieve this, we define the bigonal construction for
singular allowable covers by {\em choosing} a flat family
of smooth covers whose limit is our allowable cover, and
{\em defining} the construction to be the limit of the
construction for the nonsingular fibers. Beauville's results
imply that this is well defined, and does give a symmetric
construction: this is more or less clear except at a
singularity of Beauville type. There the problem reduces to
a local calculation whose answer, which we record in
\ref{type}. below, is visibly symmetric.

We will need a few properties of this construction (see
\cite[Section 2.3]{don})
\begin{enumerate}
\item As we said, the construction over $U$ is symmetric:
starting with $\tilde{{C'}},\dots,f'$ gives back $\Ctil,\dots,f$.
\item
\label{type}
Denote the type of $\Ctil/C$ at a point $k\in K$ by
\begin{itemize}
\item $\bto$ if $C$ is unramified over $k$ and $\Ctil$ is
ramified over exactly one point in $f^{-1}(k)$;
\item  $\btt$  if $C$ is ramified
over $k$ but $\Ctil$ is unramified over the point
$f^{-1}(k)$;
\item $\btth$ if $C$ is
unramified over $k$ and $\Ctil$ is ramified over both
branches of $C$ over $k$;
\item $\btf$ if both $C$,
$\Ctil$ have ordinary double points above it, and $\Ctil/C$
is of Beauville type there.
\end{itemize}
If $\Ctil/C$ is of type $\bto$, $\btt$, $\btth$, $\btf$ at
$k$ then $\tilde{{C'}}/C'$ is respectively of type
$\btt$, $\bto$, $\btf$,$\btth$ there.

Notice that normalization takes type $\btf$ to type $\btth$.
\item The natural $2$-$2$ correspondence between $\Ctil$
and $\tilde{{C'}}$ induces an isogeny
$\PrymC\ra\PrymCp$, whose kernel is the same as the kernel
of the natural isogeny
$\PrymC\ra\PrymC^\vee$ induced by the polarization from
$\PrymC$ to its dual abelian variety $\PrymC^{\vee}$.
In other words we get an isomorphism 
$\PrymC^\vee\sra{\sim}\PrymCp$ (cf. Pantazis \cite{pan},
at least
when $K=\PP^1$ which is all we need). As a check, let
$a$, $b$, $c$, and $d$
be the numbers of points where $\Ctil/C$ is of type
$\bto$, $\btt$, $\btth$, and $\btf$ respectively.
Then by Lemma~\ref{polar} the polarization type
for $\PrymC$ is $1^{\frac{a+2c}{2}-1}2^{\frac{b+2d}{2}-1}$.
Similarly the polarization type for
$\PrymCp$ is obtained by interchanging $a$ with $b$ and
$c$ with $d$, and this gives exactly the type dual to the
one of $\PrymC$.
\end{enumerate}

For Recillas's trigonal construction start with
$K$, $C$, $\Ctil$, $\pi$, and $f$ as before except that $f$
now has degree $3$. We get a cover $g:X\ra K$
of degree $4$ by making over the smooth unramified part
$U$, defined as before, a
construction analogous to what we previously did to get $C'$.
Namely, let $\tilde{X}/U$ represent the sheaf of sections of
$\pi:  (f\pi)^{-1}U\ra f^{-1}U$, and define $X/U$ as the
quotient of $\tilde{X}$ divided by $\iota$ (which is
defined as before). In the nonsingular case we define $X$
as the closure of $X/U$ in $X\times X\times X$, and in the
general allowable case by taking a flat limit of the construction
for smooth, unramified covers. Here we have
(\cite[Section 2.4]{don})
\begin{enumerate}
\item Over $U$ the construction is reversible: $C_{|U}$
represents the sheaf of partitions of $X_{|U}$ to two pairs of
sections, and $\tilde{{C'}}$ represents the choice of one of these
pairs.
\item Denote the type of $\Ctil/C$ at a point $k\in K$
by
\begin{itemize}
\item $\tto$ for $C$, $\Ctil$ if $C$ has exactly one simple
branch point over $K$ and $\pi$
is unramified over $f^{-1}(k)$;
\item $\ttt$ if $f$ is unramified at $k$ and $\pi$ is branched
over two of the branches of $f$ and unramified over the third;
\item $\tttp$ if two branches of
$C$ over $K$ cross normally, the third is unramified, and
moreover, if $\Ctil/C$ is of Beauville type over the double
point and unramified over the unramified branch.
\end{itemize}
Then $X$ has exactly one simple branch point at a point
$k\in K$ of type $\tto$ for $\Ctil/C$, and we denote
by $\ttoX$ the type of $X$ over $k$. Conversely, if $X$
is of type $\ttoX$ at $k$
then $\Ctil/C$ is of type $\tto$ there. If $\Ctil/C$ is
of type $\ttt$ at $k$ then $X$ has two simple branch
points over $k$, which we denote by type $\tttX$.
Here the situation is not reversible:
if $X$ is of type $\tttX$ at $k$ then $\Ctil/C$ is of
type $\tttp$ there. Notice that normalization takes type
$\tttp$ to type $\ttt$.
\item If $K\simeq\PP^1$ and $\Ctil/C$ is allowable, then
$X$ is smooth and $\Jac(X)\simeq\Prym(\Ctil/C)$. This is
due to Recillas when $\Ctil/C$ is smooth unramified, and
again limiting arguments imply this in general.
\end{enumerate}

\section{The genus $2$ case}
\label{genustwo}
Humbert's correspondence of curves of genus $2$ was studied
 by Bost and Mestre (see \cite{hum},
\cite{bome}). We shall show how to make this correspondence
via the bigonal construction, and use this to determine the
type of the isogeny.

Humbert's construction starts with a conic $C$ in $\PP^2$
with $6$ general points on it (see Remark~\ref{gen} below),
which are given as $3$
unordered pairs $\{P'_i,P_i''\}$, $i=1,2,3$. It associates
to these $3$ new unordered
pairs of points, all distinct, on $C$ as follows. Let
$\ov{P'_iP_i''}$ be the $3$ lines joining paired points,
and let  $l_k$ be the intersection of $\ov{P'_iP_i''}$ and
$\ov{P'_jP_j''}$ if $\{i,j,k\}=\{1,2,3\}$. The new
$3$ unordered pairs of points on $C$ are then the pairs of
points of tangency to $C$ from the $l_k$\/'s.

For our purposes it is more convenient to view the new points
as lying on the conic $C^*$ dual to $C$ in the dual plane
$\PP^{2*}$. A point of $\PP^{2*}$ is a line in $\PP^2$; it
is in $C^*$ if and only if this line is tangent to $C$. Let
$\phi:\PP^2\ra\PP^{2*}$ be the isomorphism defined by $C$;
namely, for $P\not\in C$ there are two tangents to $C$
through $P$, and
$\phi(P)$ is the line joining their points of tangency. For
$P\in C$, $\phi(P)$ is the tangent to $C$ at $P$. Under the
isomorphism $\phi_{|C}:C\sra{\sim} C^*$, Humbert's new pairs
go to the pairs $L'_k,L''_k$ of tangents to $C$ through $l_k$.
\begin{theorem}
Let $\pi:H\ra C$ and $\pi^*:H^*\ra C^*$ be double covers
branched over the old and new sets of points respectively. Then
there is an isogeny $\Jac(H)\ra\Jac(H^*)$ whose kernel
is a lagrangian subgroup of $\Jac(H)[2]$.
\end{theorem}
\begin{pf} Choose some $k\in\{1,2,3\}$. The set $L^* = L_k^*$ of
lines through $l_k$ is the line in $\PP^{2*}$ dual to $l_k$.
Let $f:C\ra L^*$ be the ``projection'' sending each point of
$C$ to the line joining it to $l_k$.
Dually, let $f^*:C^*\ra L = \ov{P_k'P_k''}$ be
the ``projection'' sending each tangent line of $C$ to its
intersection with $L$. Let $\psi:L\ra L^*$ be the
isomorphism sending a line through $l_k$ to its intersection
with $L$. The maps $f$ and $f^*$ have degree $2$, and hence
also $g=\psi f^*$ has degree $2$. Both coverings
$H\sra{\pi}C\sra{f}L^*$ and $H^*\sra{\pi^*}C^*\sra{g}L^*$
are unramified over the complement in $L^*$ of the six points
\begin{itemize}
\item The tangents $L'_k,L''_k$ to $C$ through $l_k$; there
$H/C$ is of type $\btt$ and $H^*/C^*$ is of type $\bto$.
\item The lines $\ov{P_i'P''_i}$ and $\ov{P_j'P''_j}$
whose intersection defines $l_k$; there both $H/C$ and
$H^*/C^*$ are of type $\btth$.
\item The lines $\ov{P'_kl_k}$ and $\ov{P''_kl_k}$;
there $H/C$ is of type $\bto$ and $H^*/C^*$ is of type
$\btt$.
\end{itemize}
It follows that if we perform the bigonal construction
on $H\ra C\ra L^*$, the two points $\ov{P_i'P''_i}$ and
$\ov{P_j'P''_j}$ of $L^*$ are of Beauville type for the
resulting cover $H'\ra C' \ra L^*$ and there are no
other singularities (see Section~\ref{bigonal}). The
preceding analysis of the ramification of $H^*/C^*$,
combined with the one for the bigonal construction
$H'/C'$ in Section~\ref{bigonal}
shows that the normalization of $H'/C'$ is isomorphic to
$H^*/C^*$. It remains to determine the kernel of the induced
isogeny on jacobians; by Pantazis's result recalled above,
it factors as
\[ \begin{array}{rcl}
\Jac H & \simeq &  \Prym(H/C) \sra{\sim} \Prym(H/C)^{\vee}
\sra{\sim} \Prym(H'/C') \\
& \sra{\nu^{\ast}} & \Prym(H^*/C^*) \simeq \Jac H^* \,.
\end{array}\]
To compute the kernel of $\nu^{\ast}$ we cannot use
Lemma~\ref{polar} directly, since $H'/C'$ is
not allowable, being ramified over two points
$x'\in g^{-1}(L'_k)$, $x''\in g^{-1}(L''_k)$. Instead
glue $x'$ to $x''$ to obtain a curve with one more double
point $C''$ and glue their inverse images in $H'$ to get
a curve $H''$, which is now an allowable cover of $C''$
($H''$ is obtained from $H$ by gluing the Weierstrass
points in pairs). We have maps of covers $H^*/C^*\ra
H'/C'\ra H''/C''$ inducing maps of Prym varieties. Applying
Lemma~\ref{polar} twice now gives that the kernel of
$\Prym(H''/C'') \ra \Prym(H/C)$ is an isotropic subgroup
isomorphic to $(\Ztwo)^2$ and that $\Prym(H''/C'')$ is isomorphic
to $\Prym(H^*/C^*)$. This implies that $\Ker\nu^{\ast}$ is as 
asserted, completing the proof of the Theorem.
\end{pf}
\begin{remark}
\label{gen} {\rm
The points $\{P'_i,P''_i\}$ are assumed general only to
guarantee that they are distinct and that the resulting
new $6$ points are also distinct (for which it suffices
that the tangents to $C$ from $l_k$ in the proof do not
touch $C$ at $P'_k$ nor at $P''_k$). In the case considered
in \cite{bome} this holds, because they assume that
$C$ and the points are real and satisfy some ordering
relations.      }
\end{remark}

\section{The hyperelliptic genus $3$ case} 
\label{hyperelliptic}
In this section we will solve our problem in the hyperelliptic
case: we will construct a correspondence between
the generic hyperelliptic curve of genus $3$ and a certain
non-generic curve of genus $3$ (which is not hyperelliptic).
Let $H$ be a hyperelliptic curve of genus $3$ and
let $\pi_1:H\ra\PP^1$ be the hyperelliptic double cover.
Choose a grouping in pairs of the $8$ branch points
$w_1,\dots,w_8\in\PP^1$ of $\pi_1$. We claim that there
exists a map $g_1:\PP^1\ra\PP^1$, of degree $3$, which
identifies paired points. This can be seen in several
ways. Firstly, let $T$ be the curve obtained from $\PP^1$
by identifying paired points to ordinary double points.
We think of $T$ as a curve of genus $4$ and take its
canonical embedding to $\PP^3$. As in the nonsingular
case, the canonical map is well behaved, and in particular the
canonical image of $T$ lies on a unique, generically
nonsingular quadric by the Riemann-Roch theorem. Projecting
via either of the two ruling of this quadric will give the
desired map $g_1$. Notice that by its construction $g_1$
factors as $\PP^1\sra{\nu} T \sra{g} \PP^1$, where $\nu$ is
a normalization map.

Another way to get $g_1$ is to embed $\PP^1$ in $\PP^3$
as a rational normal curve. We look for a projection
from $\PP^3$ to $\PP^1$ which identifies paired points.
The center of this projection is a line $L$ which
must meet the $4$ lines joining the pairs. The
grassmanian $G(1,\PP^3)$ of lines in $\PP^3$ is
naturally a quadric in $\PP^5$ and the condition to meet
a line is a linear condition. We see again that there is
always at least one such $L$, and generically two.

We now perform the trigonal construction. This gives a map
of degree $4$ $f:C\ra\PP^1$ sitting in a diagram
\begin{equation}
\label{hyp1}
\begin{array}{rcl}
&&H \\
&&\da \! \spione \\
C&&\PP^1\\
&\sf\!\!\sea\;\;\swa\!\!\sgone&\\
&\PP^1&   \end{array}
\end{equation}
Let $w_{12},\dots,w_{78}$ be the $4$ images of the $w_i$\/'s
under $g_1$, with the indices indicating the grouping.  By
the Riemann-Hurwitz formula there are generically $4$ points
$a_1,\dots,a_4$ in $\PP^1$ over which $g_1$ is branched, with
a simple branch point over each. Hence $H/\PP^1$ is of type
$\tto$ at each $a_i$ and of type $\ttt$ at each $w_{2i-1,2i}$.
>From the properties of the trigonal construction we get
$2-2g(C)=8-8-4$, so that $C$ has genus $3$.
The trigonal construction gives a birational correspondence between
\begin{itemize}
\item The moduli of the data
$(H \sra{\pi_1}\PP^1 \sra{g_1} \PP^1)$ with
$4$ points of type $\tto$ and $4$ points of type $\ttt$.
\item A component of the Hurwitz scheme parametrizing
$4$-sheeted covers $f:C\ra \PP^1$ with $4$ simple branch
points and $4$ double branch points.
\end{itemize}

Each of these moduli spaces is $5$ dimensional. (Another
component of this Hurwitz scheme parametrizes bielliptics,
namely maps $f:C\ra \PP^1$ which factor through a double cover
$E \ra \PP^1$ where $E$ is elliptic. Curves in this latter
component are taken by the trigonal construction to towers
$H\ra \overline{H} \ra \PP^1$ where $\overline{H}=A \cup B$
is reducible, with $A$, $B$ of degrees $1$, $2$ respectively
over $\PP^1$. We shall not need this component in what follows.)

The key point for us is that the trigonal construction induces
an isogeny $\Jac(C)\ra\Jac(H)$ whose kernel is lagrangian
in $\Jac(C)[2]$. More precisely we have the following
\begin{proposition}
Let $\PP^1\sra{\nu}T\sra{g}\PP^1$ be as before, and let $\Ttil$
be the curve obtained by identifying the Weierstrass points in
$H$ to ordinary double points with the same grouping as the one
we chose to get $T$. Then
\begin{enumerate}
\item Diagram (\ref{hyp1}) extends to
\[ \begin{array}{rcccl}
&&\Ttil&\sla{\nutil}&H\\
&&\spi\!\da\;&&\;\da\!\snu\spi\\
C&&T&\sla{\nu}&\PP^1\\
&\sf\!\!\sea\;\;\swa\!\!\sg&\\
&\PP^1&   \end{array}\,. \]
Here $\nutil:H\ra\Ttil$ is the normalization map, and we view
$\nu\pi:=\pi_1:H\ra \PP^1$ as the normalization of
$\pi:\Ttil\ra T$.
\item $\nutil$ induces an isogeny of polarized abelian varieties
$\nutil^*:\Prym(\Ttil/T)\ra\Jac(H)$ whose kernel is
lagrangian in $\Prym(\Ttil/T)[2]$.
\item Let $\phi:\Jac(C)\ra\Jac(H)$ be the isogeny obtained
by composing $\nutil^*$ with the isomorphism
$\Jac(C)\simeq\Prym(\Ttil/T)$.
Then the kernel of $\phi$ is lagrangian in $\Jac(C)[2]$, and the 
kernel of the dual isogeny $\phi^*:\Jac(H)\ra\Jac(C)$ is the 
lagrangian subgroup of $\Jac(H)[2]$ generated by the differences of 
identified Weierstrass points.
\end{enumerate}
\end{proposition}
\begin{pf}
Part 1. holds because $\Ttil/T$ is allowable. The pairs
of points of $H$ identified by $\nu$ lie over points of
type $\tttp$ for $T$, $\Ttil$. Hence they are branch points
for $\pi$, namely Weierstrass points. The rest follows from
Lemma~\ref{polar}.
\end{pf}

\section{The generic genus $3$ case}
\label{genusthree}
Let $C$ be a generic curve of genus $3$. In this section we shall
give a construction of a curve $C'$ of genus $3$ and an isomorphism
$\Jac (C)/L \simeq \Jac (C')$ where $L$ is a lagrangian subgroup of
$\Jac (C)[2]$. Let $f:C\ra\PP^1$
be a map of degree $4$, and let $b_1$, $b_2$ be points in
$\PP^1$ such that $f$ has two simple branch points over each
$b_i$. It is easy to show such $f$, $b_1,b_2$ exist, and in fact we will  
parametrize the space of such $f$\/'s in the end of this section.

We perform the trigonal construction on $f$. This gives curves
$T$, $\Ttil$ and maps $g:T\ra\PP^1$ and $\pi:\Ttil\ra T$, with
$\deg g = 3$ and $\deg \pi=2$. Let $\nutil:\nu \Ttil\ra\Ttil$
and $\nu:\nu T\ra T$ be normalization maps and let
$\nu\pi:\nu\Ttil\ra\nu T$ be the map induced by $\pi$. The
properties of the trigonal construction show the following.
Firstly, $T$ and $\Ttil$ have each two ordinary double points,
one over each $b_i$, and no other singularities. Next, the
map $g\nu:\nu T\ra\PP^1$ has exactly $8$ branch points, all
simple, one over each $a_i$. It follows that the genus
$g(\nu T)$ is $2$ and therefore the arithmetic genus $g(T)$
is $4$. The map $\nu\pi$ has exactly $4$ ramification points
$P_i$, $Q_i$, two over each $b_i$ for $i = 1,2$, and hence
$g(\nu\Ttil)=5$ and $g(\Ttil)=7$.

Since $\nu T$ has genus $2$, it is hyperelliptic. Let
$h:\nu T\ra \PP^1$ be the hyperelliptic double cover, and
let $w_1,\dots,w_6\in\PP^1$ be the branch points of $h$.
The bigonal
construction gives curves and maps of degree $2$
$\nu\Ttil'\sra{\nu\pi'}\nu T'\sra{h'}\PP^1$. The points
in $\PP^1$ over which $h\nu\pi$ is not \'etale are the $6$ 
$w_i$\/'s, which are of type $\btt$ for $\nu T$ and $\nu\Ttil$, 
and the $4$ points $h(P_i)$\/, $h(Q_i)$\/,  $i=1,2$, which
are of type $\bto$. The types get reversed for $\nu T'$ and
$\nu\Ttil'$, and in particular $\nu T'$ is ramified exactly
over the $h(P_i)$\/'s and the $h(Q_i)$\/'s. It follows that  
$g(\nu T')=1$. We also see that $\nu\pi'$ has $6$ branch points,
say $w'_1,\dots,w'_6$, one over each of the $w_i$\/'s, and
hence $g(\nu\Ttil')=4$. The curves $\nu T'$ and $\nu\Ttil'$
are nonsingular.

Choose a grouping of the $w_i$\/'s in $3$ pairs. Identify
the corresponding $w'_i$\/'s in $\nu T'$ to get a curve
$T'$ with $3$ ordinary double points, say $w'_{12}$, $w'_{34}$,
$w'_{56}$, the indices indicating the groupings. $T'$ has 
arithmetic genus $4$. Likewise identify the corresponding
points above the $w'_i$\/'s on $\nu\Ttil'$ to obtain a curve
$\Ttil'$ with $3$ ordinary double points and arithmetic
genus $7$.

As in the nonsingular case, the canonical embedding sends $T'$ 
to $\PP^3$ and the image sits on a unique, generically smooth
quadric. Choosing one of the two rulings of this quadric gives
a map $g':T'\ra\PP^1$. This map is of degree $3$, because the
canonical curve is a curve of type $(3,3)$ on the quadric. The
map $g'\nu':\nu T'\ra\PP^1$ is ramified over $n=6$ points, since
$2-2g(\nu T')=0=3(2-2g(\PP^1))-n$. Over these the pair $\nu T'$,
$\nu\Ttil'$ is of type  $\tto$. There are also $3$ points of
type $\ttt$, the images under $g'$ of the identified pairs $w'_{12}$,
$w'_{34}$, $w'_{56}$.  The
trigonal construction performed on
$\nu\Ttil'\slra{\nu\pi'}\nu T'\slra{g'\nu'}\PP^1$ gives a
curve $C'$ and a map $f':C'\ra\PP^1$ of degree $4$. We
readily see it has genus $3$. The following diagram summarizes
the procedure:

\[  \begin{array}{rcccccccccl}
&&{}_7\Ttil&\sla{\snutil}&{}_5\nu\Ttil&&{}_4\nu\Ttil'&
\sra{\snutil'}& {}_7\Ttil'&&\\
&&\spi\!\da&&\da\snu\spi&&\snu\spi'\!\da&&\da\spi'&&\\
{}_3 C&&{}_4 T&\sla{\nu}&{}_2\nu T&&{}_1\nu T'&\sra{\nu'}&
{}_4 T'&&{}_3 C'\\ 
&\sf\!\sea\;\swa\!\sg&&&& \sh\!\sea\;\swa\!\sh'&&&& 
\sg'\!\sea\;\swa\!\sf'&\\
&\PP^1&&&&\PP^1&&&&\PP^1&
\end{array}\,.    \]

Before stating our main result we need to discuss the choices
made in the construction. Writing $f^{-1}(b_i) = 2(P_i+Q_i)$,
we obtain a point of order $2$
\[ \alpha = \alpha(f) = P_1 + Q_1 - P_2 - Q_2.  \]
in $\Jac(C)$. The trigonal isomorphism 
$\Jac(C) \simeq \Prym(\Ttil/T)$ maps $\alpha$ to the 
difference $L(b_1) -L(b_2)$ (defined in the discussion 
preceeding Lemma \ref{polar}), which is the nontrivial
element in $\Ker\nu^\ast$. Moreover the only choice 
made other than $f$ is the grouping of $w_1,\dots,w_6$ 
under $\nu'$. The differences of the corresponding 
paired points in $\nu T$ are the nonzero elements of a 
lagrangian subgroup $L_0$ of $\Jac(\nu T)[2]$.

Now observe that the pullback to $\nu\Ttil$ by $\nu\pi$ of a
line bundle of order $2$ on $\nu T$ is in the kernel of the
norm map to $\nu T$. This gives a symplectic embedding
\[ \iota:\Jac(\nu T)[2] \hookrightarrow \Prym(\nu\Ttil/\nu T)[2]\,. \]
The image $\Ima(\iota)$ of $\iota$ can be described in two ways.
On the one hand, it is the kernel of the polarization map
$\Prym(\nu\Ttil/\nu T) \ra  \Prym(\nu\Ttil/\nu T)^\vee$.
(Observe that $\Prym(\nu\Ttil/\nu T)$ has a polarization
of type $221$ by Lemma~\ref{polar}, whose kernel is then
isomorphic to $(\ZZ/2\ZZ)^4$.) On the other hand, let
$\alpha^\perp$ denote the orthogonal complement to (the image
of) $\alpha$  in $\Prym(\Ttil/T)$ for the Weil
pairing $w_2$. Then $\Ima(\iota)$ is
also the pullback of $\alpha^\perp$ by the normalization 
map. Indeed, for $u\in \Jac(\nu T)[2]$ we have 
$\iota(u)\in\alpha^\perp$ because
\[ w_2(\iota(u),\alpha) = w_2(u,\Nm_{\Ttil/T}(\alpha)) = 0\,, \]
and as both groups have cardinality $16$ they coincide. Hence
this image is isomorphic to
$\alpha^\perp / \langle\alpha\rangle$. In particular, the 
inverse image $L$ of $L_0$ in $\Jac(C)$ is a lagrangian 
subgroup of $\Jac(C)[2]$ containing $\alpha$.

Conversely, let $L\subset \Jac(C)[2]$ be a lagrangian
subgroup. We will say that the choices $f$, $\nu'$ made
in the course of the construction are compatible with $L$
if $\alpha = \alpha(f)$ is in $L$ and $\nu'$ corresponds 
to $L/\langle \alpha \rangle$ as above.
We can now formulate our main theorem, to which we shall 
give two proofs:
\begin{theorem}
\label{main}
Let $C'$ be the result of the construction applied to a curve 
$C$ of genus $3$ compatibly with a lagrangian subgroup
$L\subset \Jac(C)[2]$. Then there is an induced isomorphism 
$Jac(C)/L \sra{\sim} \Jac(C')$. In particular $C'$ is 
independent of the (compatible) choices made in the construction.
\end{theorem}
\begin{pf}

The construction induces isogenies whose degrees are marked
below:
\[
\begin{array}{rcl}
 \Jac(C) & \simeq & \Prym(\Ttil/T) \sbth{\nu^\ast}{\lra}{2}
\Prym(\nu\Ttil/\nu T) \sbth{\delta}{\lla}{4}
{\Prym(\nu \Ttil'/\nu T')}\\
&\sbth{{\nu'}^{\ast}}{\lla}{4}& \Prym(\Ttil'/T')
\simeq  \Jac(C')\,.
\end{array}
\]
Here the middle step $\delta$ is identified with the 
polarization map from an abelian variety of polarization
type $211$ to its dual. As before we identify $\nu^\ast$ 
with the quotient by $\alpha$, so to construct our isomorphism
$\Jac(C)/L \simeq \Jac(C')$ it would suffice to produce
a natural map
$\epsilon: \Prym(\nu\Ttil/\nu T) \ra \Prym(\Ttil{}'/ T')$
whose kernel is the subgroup $L/\langle \alpha \rangle$ of
$\Prym(\nu\Ttil/\nu T)$.

One way to do this is to define $\epsilon$ as the dual map
of $\nu'{}^\ast$, using $\delta$ to identify the dual of
$\Prym(\nu\Ttil'/\nu T')$ with $\Prym(\nu\Ttil/\nu T)$, and
using the principal polarization on $\Prym(\Ttil'/T')$ to
view it as its own dual. Tracing through the definitions
one verifies that $\Ker(\epsilon)$ is indeed
$L/\langle \alpha \rangle$ as asserted.

An alternative, and more geometric approach, is to show that 
the hyperelliptic case treated in Section~\ref{hyperelliptic}
is a specialization of our present general construction. In 
fact the hyperelliptic case is obtained when the $4$-sheeted 
cover $f: C\ra \PP^1$ happens to have $4$, rather than the
generic $2$, double branch points. We shall see that this
determines a preferred gluing $\nu'$. In going to this
special case we have to note that the limits of the curves
$\nu T$, $\nu \Ttil$ (which we continue to denote with the
same symbols) are no longer non-singular: they are now
only partial normalizations of $T$, $\Ttil$, and the map
$\nu \Ttil \ra \nu T$ now has $2$ points of Beauville type,
at the singularities which were not normalized. The full
normalizations, say $\nu \nu T$ and $\nu \nu \Ttil$, now have
genera $0$ and $3$ respectively, and the resulting diagram
\[\begin{array}{rcl}
&& \nu\nu\Ttil \\
&& \spi\!\da \\
C && \nu\nu T \\
&\sf\!\!\sea\;\!\swa\!\!\sg' &\\
& \PP^1 &
\end{array} \]
clearly coincides with diagram~(\ref{hyp1}).

Reagardless of the singularities of the intermediate curves,
we will see that each of the abelian varieties in the
diagram specializes to an abelian variety. In particular,
the limit of $\Prym(\nu\Ttil/\nu T)$ is, by Lemma~\ref{polar},
a $4$\/-sheeted cover of
$\Prym(\nu\nu\Ttil/\nu\nu T) \simeq \Jac(H)$, whose kernel is
$L/\langle \alpha \rangle$. Below we will also
identify the limit of $\Prym(\Ttil{}'/T')$ with $\Jac(H)$.
This will produce the desired map $\epsilon$ in this special
case, and hence in general.

For this, we note that the bigonal data
$\nu \Ttil \ra \nu T \ra \PP^1$ has $4$, $2$, $0$ and $2$
points of types $\bto$, $\btt$, $\btth$, and $\btf$
respectively, which turn into points of types $\btt$, $\bto$,
$\btf$ and $\btth$, respectively, for
$\nu \Ttilp \ra \nu T' \ra \PP^1$.
To obtain $\Ttilp \ra T'$ we need to pair the $6$ ramification
points of $\nu \Ttilp \ra \nu T'$. There are $15$ ways to do
this, of which one is distinguished: each pair of Beauville
branches gets paired, as do the remaining two ramification
points. Let $\Tph$ be the intermediate object, obtained by 
gluing only the Beauville branches but not the remaining
pair. It is a singular hyperelliptic curve of genus $3$, and
is a partial normalization of $T'$ at the double point $p$.
Let $\Ttilph$ be the corresponding $1$-point partial 
normalization of $\Ttilp$, of arithmetic genus $6$.
 
To continue our construction, we need to identify the two
$g_3^1$'s on $T'$: these two turn out to coincide, and the
unique $g_3^1$ is in fact given by the $g_2^1$ on $\Tph$ 
plus a base point at the double point $p$. To see this we
examine what happens to our general construction of the 
$g_3^1$ in this case. The unique quadric surface through
the canonical model of $T'$ is now a quadric cone, with
vertex at (the image of) $p$, because projection from
$p$ gives the canonical image of the hyperelliptic $\Tph$,
which is the double cover of a conic. Therefore the two
rulings, hence the two $g_3^1$'s, coincide and have a
base point at $p$, as asserted.

At this point we need to turn the $g_3^1$ into a morphism, 
which requires us to blow up the point $p$. This
results in a reducible trigonal curve $\Tpt := \Tph \cup P$,
where $P$ is a copy of $\PP^1$ intersecting $\Tph$ in the two
inverse images $p_1$, $p_2$ of $p$ in $\Tph$. The trigonal
map has degrees $2$ and $1$ respectively on the two 
components $\Tph$ and $P$. This curve is indeed a flat limit, 
in the family of triple covers of $\PP^1$, of the trigonal
curves encountered in the non-hyperelliptic situation. The
corresponding double cover $\Ttilpt \ra \Tpt$ is of Beauville
type at all $4$ of the singular points (the two singularities
of $\Tph$ plus $p_1$, $p_2$). Here
$\Ttilpt = \Ttilph \cup \Ptil$, where $\Ptil$ is another copy
of $\PP^1$, double cover of $P$ branched at the points glued
to $p_1$ and $p_2$.

Now that we have identified the trigonal data, we can complete
the construction. By example~2.10(iii) of \cite{don}, or by 
inspection, we see that the result $C'$ of applying the 
trigonal construction to the reducible trigonal data
$(\Ttilph \cup \Ptil) \ra (\Tph \cup P) \ra \PP^1$ is the
$4$-sheeted cover of $\PP^1$ obtained by applying the bigonal 
construction to $\Ttilph \ra \Tph \ra \PP^1$. But since the 
bigonal construction is reversible this is nothing but the
hyperelliptic curve $H= \nu\nu \Ttil$ which resulted from 
the construction of Section~\ref{hyperelliptic}, as claimed.

The degeneration just described involves a flat family of
abelian varieties, so the polarization type and the type of 
the kernel of the isogeny  on jacobians remain constant. From
the hyperelliptic case we now see that $L$ is the kernel of
our isogeny in the general case. By Torelli's theorem, $C'$
is determined by its polarized jacobian, which is $\Jac(C)/L$.
Hence $C'$ is indeed independent of the choices (compatible 
with $L$) made during the construction. This concludes the 
proof of Theorem~\ref{main}.
\end{pf}

We now make some further comments on the
choices we made in the course of the construction. 
Starting on the left, we fix the curve C, the Lagrangian 
subgroup L and an element $\alpha \in L$. Our 
$g^1_4$\/'s $f: C \ra \PP^1$ with $\alpha =\alpha(f)$ 
are determined by a divisor class in the intersection
\[Z = Z_\alpha =
\Theta_C \cap(\alpha + \Theta_C) \subset \Pic^2(C)\,.\]
More accurately $Z$ parametrizes the family of $g^1_4$'s
(with the specified $\alpha$), together with a marking of
the two singular points $b_1$, $b_2$. Interchanging these
two points gives an involution $i$ of $Z$ induced by the
involution $x \mapsto x+ \alpha$ of $\Pic^2(C)$, and it
is the quotient of $Z$ by $i$ which parametrizes the
$g^1_4$'s alone. $Z$ is also invariant under the involution
$j: x \mapsto K_C- x$ (where $K_C$ is the canonical class)
and $i$ and $j$ commute. In addition, since $\Theta_C$ is 
an ample divisor, $Z$ is connected.
Counting fixed points shows that the respective quotients of 
$Z$ by $i$, $j$, $k=ij$ have genera $4$, $1$, $4$, and that the
common quotient $\ov{\ov{Z}}:=Z/\langle i,j\rangle$ has genus $1$.




These quotients clearly have the following interpretations as 
parameter spaces:

\begin{enumerate}
\item $Z$ parametrizes the $g^1_4$\/'s $f:C\ra \PP^1$ (equivalently,
via the trigonal construction, 
towers of double covers $\Ttil \ra T \ra \PP^1$ of the indicated type), 
{\rm with} a choice of a double branch point $b_1$.
\item $Z/i$ parametrizes the $g^1_4$\/'s $f:C\ra \PP^1$ (equivalently,
towers of double covers $\Ttil \ra T\ra\PP^1$ of the indicated type).
\item $Z/j$ parametrizes the double covers $\Ttil\ra T$ of the indicated 
type together with a singular point of $T$.
\item $\ov{\ov{Z}}$ parametrizes the double covers $\Ttil\ra T$ of the 
indicated type, hence it also parametrizes their normalizations, as well 
as the maps $\nu\pi': \nu\Ttil' \ra  \nu T'$.
\end{enumerate}

We will now discuss what choices we make when we perform the
construction in reverse order, and how the choices from the
two directions are related. 

Starting with the genus $3$ curve $C'$, we now assume given a lagrangian 
subgroup $L'$\/($=\Jac(C)[2]/L)$ of $\Jac(C')[2]$, and a subgroup
$G \subset L'$ of order 4 (which corresponds to 
$\alpha^\perp/L$). A marking 
of the three double branch points $b'_i$ of $f'$  is equivalent to a 
choice of a basis  $\beta' = P'_2 + Q'_2 - P'_1 -Q'_1$ and
$\gamma' = P'_3 + Q'_3 - P'_1 -Q'_1$ of $G$. Let 
$\Theta' \subset \Pic^2(C')$ be the theta divisor of $C'$, and for 
a class  $u\in\Jac(C')$ let $\Theta'_u$ denote the translation of 
$\Theta'$ by $u$. Consider a line bundle $L$ in the intersection
$S = \Theta \cap \Theta_{\beta'} \cap \Theta_{\gamma'}$. 
Since the canonical bundle is the only degree 4 bundle on $C'$ 
with $h^0>2$, there are only two possibilities: either $L^{\otimes 2}$ gives a 
$g^1_4$\/ $f':C'\ra \PP^1$ with three marked
double branch points $b'_i\in \PP^1$, $i=1,\dots,3$, or else $L$ 
must be a theta characteristic on $C'$. We claim 
the following:\\

\noindent (1) $S$ consists of six points \\

\noindent (2) $S$ is closed under $v\ra K_{C'} - v$.\\

\noindent (3) Four of the points of $S$ are theta characteristics, and two are not.\\

\begin{pf}
(1) holds because $6=g!$.  For (2), suppose that $f',f'':C' \ra \PP^1$
correspond to $2v$, $2K_{C'} - 2v$ respectively. Then for each double 
ramification point $P'_i$, $Q'_i$ of $f'$ we get a unique double 
ramification point $P''_i$, $Q''_i$ for $f''$ by imposing the condition
$P''_i + Q''_i + P'_i + Q'_i = K_{C'}$. 

For (3) , one checks that there is a unique coset $G'$ of 
$G$ in the set of odd theta characteristics
on $C'$: indeed, in coordinates we may take the set of theta characteristics 
to be $V = (\ZZ/2\ZZ)^6$ with coordinates $x_1,\dots,x_6$, and we may suppose that  $h^0(C',O_{C'}(x))$ mod 2 for $x = (x_1,\dots,x_6)$ is given by 
$q(x) = x_1x_2 + x_3x_4 + x_5x_6$. Also we may simultaneously identify 
$\Jac(C')[2]$ with $V$, with the Weil pairing given by
$w_2(x,y) = q(x+y)-q(x) -q(y)$. Without loss of generality we can also
take $\beta' = e_1$ and $\gamma' = e_3$, with $e_i$ the standard $i$\/th
unit vector. Then  $G' = \{(a,0,b,0,1,1)\}$. 

Part (3) is now clear:   $G'$ is contained in $S$, and no other theta characteristics 
appear in $S$. This establishes our claim. 
\end{pf}

We can now describe all the choices made when we start from the right side. 
Our data $C',L',G$ determines a complementary pair of maps $f',f"$. These 
determine the data 
$\Ttil'\sra{\pi'} T'$ uniquely (the two resulting maps $g',g"$ are the 
usual two $g^1_3$'s on the genus 4 curve $T'$). The normalization 
$\nu\pi': \nu\Ttil' \ra  \nu T'$ is therefore also uniquely
determined. So the {\em only} choice made is that of $h'$, given by 
an arbitrary point of 
$ Pic^2(\nu T') \approx  \nu T'$. Comparing with what we found starting 
from the left, we discover that $\nu T'$ is precisely identified with the double 
quotient $\ov{\ov{Z}}$.



\section{The case of genus$\geq 4$.}
One might try to generalize our
construction to higher genus by finding, for a generic curve
$C$ of genus $g$, a correspondence with another generic curve
$C'$ of genus $g$ such that $\Jac(C')\simeq\Jac(C)/K$, with
$K$ a lagrangian subgroup of $\Jac(C)[2]$. We shall show
that this is not possible.
\begin{theorem}
Let $K$ be a lagrangian subgroup in $\Jac(C)[p]$, where $C$ is a
generic curve of genus $g\geq 4$ and $p$ is a prime. Then
$\Jac(C)/K$ with its induced principal polarization is not a jacobian.
\end{theorem}
\begin{pf}
Let $\cT$, $\cS$, $\cM$ and $\cA$ denote respectively the 
Teichm\"uller space, the Siegel space, the moduli space of
curves and the moduli space of principally polarized 
abelian varieties, all of genus $g$. The mapping class group
$M=M(g)$ acts on $\cT$ with quotient $\cM$ and the modular 
group $\Gamma=\Sp(2g,\ZZ)$ acts on $\cS$ with quotient $\cA$. 
Moreover $\Gamma$ is naturally a quotient of $M$, because 
$M$ acts on symplectic bases for $H_1(C,\ZZ)$ through its
action on $\pi_1(C)$, and the period map $\tau:\cT\ra\cS$
is $M$-equivariant for these actions. Passing to the
quotient, we get Torelli's map $\taub:\cM\ra\cA$, which is
injective (Torelli's theorem) and exhibits $\cM$ as a locally closed
subvariety of $\cA$. Since $\cT$ is irreducible it follows
that the {\em Torelli space} $\cTb=\tau(\cT)$ is a locally closed
irreducible analytic subvariety of $\cS$.

\[ \begin{array}{rcl}
   \cT & \slra{\tau} & \cS \supset \ov{{\cT}}=\tau (\cT) \\
   M \da && \da \Gamma \\
   \cM & \slra{\ov{{\tau}}} & \cA
   \end{array}  \]

Let $W$ be the finite cover
of $\cS$ obtained by taking over each marked abelian
variety $A$ the lagrangian subgroups of $A[p]$. Since $W$ is
unramified over the contractible space $\cS$, it is in fact a
union of copies of $\cS$. Our generic isogeny
$\Jac(C)\ra\Jac(C)/K$ translates to the following data. The
curve $C$ lives over an open subset of $\cM$, hence of $\cTb$.
The subgroup $K$ corresponds to a sheet of $W_{|\cTb}$. 
Therefore our isogeny extends to the quotient map by
the subgroup, still denoted $K$, corresponding to the
``same'' sheet over all of $\cS$. Now recall that $\cS$
is the  space of symmetric
$g\times g$ complex matrices $\Omega$ with positive imaginary
part, and the abelian variety over $\Omega$ is
$A_\Omega = \CC^g/(\ZZ^g+\Omega\ZZ^g)$. Since  monodromy
(i.e. $\Gamma$) acts transitively on the lagrangian subgroups of
$A_\Omega$, we may
take $K = (\frac{1}{p}\ZZ/\ZZ)^g$ for convenience. Then
$A_\Omega/K \simeq A_{p\Omega} = A_{s\Omega}$, with
\[ s = \matr{pI_{g\times g}}{0}{0}{I_{g\times g}}
                                 \in\Sp(2g,\RR)\,.\] 

If $\Jac(C)/K$, with its principal polarization, were a jacobian,
it would follow that the Torelli locus $\cTb$ was invariant
under the subgroup $\Delta$ of $\Sp(2g,\RR)$ generated by
$\Gamma$ and by $s$. We claim that $\Delta$ is dense
in $\Sp(2g,\RR)$. Indeed, consider the subgroup $N_+$ of $\Sp(2g, \RR)$
consisting of the matrices
$n(x) = \matr{I_{g\times g}}{x}{0}{I_{g\times g}}$, where
$x$ runs over the real symmetric $g\times g$ matrices. Then
$\Delta$ contains $s^i n(x) s^{-i}=n(p^{-i}x)$ for all integral
symmetric matrices $x$ and integers  $i$. These are dense in $N_+$,
and $\Delta$ likewise contains a
dense subgroup of $N_-={}^t N_+$. It is
well-known (and easy) that $N_+$ and $N_-$ generate $\Sp(2g,\RR)$, so
$\Delta$ is indeed dense in $\Sp(2g, \RR)$.

Therefore, under our assumption, $\cTb$ would be dense in 
$\cS$ (in the complex topology),
so that $\cM$ would be dense in $\cA$. This is a contradiction
when $g>3$, because for dimension reasons
$\cM$ is not
dense in $A$ even for the Zariski topology then,
and the theorem follows.
\end{pf}

\newlength{\largeur}
\settowidth{\largeur}{Hebrew university of Jerusalem}
\parbox[t]{\largeur}{Mathematics department\\
University of Pennsylvania\\
Philadelphia, PA 19104-6395\\
USA\\
donagi@math.upenn.edu}
\hfill
\parbox[t]{\largeur}{Institute of mathematics\\
Hebrew university of Jerusalem \\
Givat Ram 91904\\
Israel\\
rlivne@math.huji.ac.il}
\end{document}